\documentclass[11pt,a4paper]{article}
\pdfoutput=1
\usepackage{amssymb}  
\usepackage{amsthm}
\usepackage{amsmath}
\usepackage{amsfonts}
\usepackage{moreverb}
\usepackage{graphicx}
\usepackage{float}
\usepackage{enumitem}
\usepackage{caption}
\usepackage{subcaption}
\usepackage{authblk}
\title{Fault Detection and Identification - a Filter Investigation}
\author[1]{Rudaba Khan}
\author[2]{Paul Williams}
\author[2]{Paul Riseborough}
\author[1]{Asha Rao}
\author[1]{Robin Hill}
\affil[1]{Department of Mathematics and Geospatial Science, RMIT University, Melbourne, Australia.}
\affil[2]{BAE SYSTEMS Australia, Melbourne Australia.}

\begin{document}

\maketitle

\begin{abstract}
This paper develops a new active fault tolerant control system based on the concept of analytical redundancy.  The novel design consists of an observation filter based fault detection and identification system integrated with a nonlinear model predictive controller.  A number of observation filters were designed, integrated with the nonlinear controller and tested before reaching the final design which comprises an unscented Kalman filter for fault detection and identification together with a nonlinear model predictive controller to form an active fault tolerant control system design.
\end{abstract}


\maketitle

\section{Introduction}
By nature a fault causes a system to behave abnormally but does not cause the system to shut down \cite{ducard2009fault}, however if left unattended a fault may lead to a system failure, hence the importance of a fault tolerant control (FTC) system.  While faults can cause instability in a system, the integration of an FTC scheme significantly increases the ability of the system to maintain overall stability in the presence of a fault \cite{zhang2008bibliographical}.  In this paper we develop a new active FTC system using an unscented Kalman filter (UKF) for fault detection and identification (FDI) where parameter updates made by the UKF are sent to the nonlinear model predictive control (NMPC) based control system.  To test the FTC system the design is applied to the control of a 2D robot model.\\       
  
There are numerous solutions to the FTC design problem with the FTC research community classifying these designs as either \textit{active} or \textit{passive}.  Passive FTC systems are based on fixed controllers and are designed to be resilient against known faults.  They are designed using robust control techniques for worst case scenarios.  Active FTC systems ``actively'' seek the fault and try to gather as much information about it as possible to help the controller overcome any consequential instabilities and are also known as self-repairing, self-designing or fault detection, identification (diagnosis) and accommodation schemes.  Active schemes are made up of an FDI component and are based on controller redesign, or selection/mixing of pre-designed controllers.\\

The main difficulty with active FTC is on-line reconfiguration which requires detailed information about changes in system parameters.  Hence the main role of the FDI subsystem is the gathering of information on parameter changes to assist in controller reconfiguration.  FDI is a key component in an active FTC system and is the most difficult aspect of FTC \cite{camacho2010fault}.  In the early days most research on FDI was done independently of the controller design and no combined design existed \cite{patton1997fault}.  Recently there has been some research on the integration of FDI and FTC; however much remains to be done.\\

Techniques commonly used for FDI include artifical intelligence based FDI schemes \cite{cork2005fault}, \cite{lin2007failure}, \cite{perhinschi2010integrated}, and multiple model based methods \cite{boskovic2001line} \cite{ducard2009fault} where different models are used to describe the dynamics of the system for different operating regimes.  Sliding mode observers for fault detection have also been used because of their strong robustness to a particular class of uncertainty \cite{edwards2010sliding}.  Other methods include analytical FDI methods that are based on estimating the fault through matrix algebra \cite{Yee2002}.  Observer based FDI remains the most common approach studied in recent literature where an extended state observer is used to estimate faults.\\

Quite often in the design of an FTC system two controllers are utilised; the main (or nominal) controller is designed for the faultless case with the second controller being a compensator designed to handle the faulty case \cite{wang2009extended}.  In this paper a number of observer based FDI schemes are investigated and integrated with the NMPC (nonlinear model predictive control) controller design of \cite{RKPaper1} to develop a full active FTC system.  In our design only one controller needs to be designed capable of handling the faultless as well as the fault cases.  The 2D robot model from \cite{RKPaper1} is used as a test bed to assist in the development of the FTC system design.\\

This paper is organised as follows; section \ref{section:chap4_FDI_TandI} presents an in-depth look at the FDI methods chosen for investigation, namely the EKF (extended Kalman filter), the UKF (unscented Kalman filter) and the IMM (integrated multiple models).  An EKF filter, UKF filter, EKF based IMM filter and a UKF based IMM filter are each designed in section \ref{section:chap4_ProbForm} for the purposes of FDI using the 2D robot model of \cite{RKPaper1}.  Four different active fault tolerant control systems using NMPC (nonlinear model predictive control) as the controller design are formulated and implemented in MATLAB for the 2D robot model.  Each of the four different active FTC systems developed are then tested under different conditions in section \ref{section:chap4_NRandA}.  The tests are designed to evaluate the performance of the filters as well as the interaction of the filter and controller designs.  Section \ref{section:chap4_NRandA} provides a detailed analysis of the test results and concludes with a summary of findings.  Based on the findings the UKF filter was chosen as part of the final design of the active FTC system.  Section \ref{subsec:LMPC} concludes with a comparison between linear MPC and nonlinear MPC as the controller showing that an active FTC system with a nonlinear MPC controller has better overall performance in the event of a failure.  Conclusions are given in section \ref{section:chap4_conclusion}.

\section{Fault Detection Techniques Selected for Implementation: Theoretical Description}\label{section:chap4_FDI_TandI}
The fault detection schemes considered here are all based on filtering techniques, namely the EKF, the UKF, the EKF based IMM and the UKF based IMM.  These filters are used to sequentially estimate the state of a dynamic system using a sequence of noisy measurements made on the system.  The state estimates are then utilised to aid in fault detection and control reconfiguration.  A general overview and key mathematical concepts are provided for each method below.

\subsection{Extended Kalman Filter (EKF)}\label{subsection:chap4_EKF_equations}
The EKF is an extension of the well known Kalman filter.  One of the drawbacks of the Kalman filter is that it does not provide good estimations for nonlinear systems \cite{ristic2004beyond}.  The EKF approximates (or linearises) the nonlinear functions in the state dynamic and measurement models.  There are two main stages during an EKF (and the general Kalman filter) cycle: predict and update.  During the prediction stage the filter states and covariances are predicted forward one time step as are the measurement predictions.  During the update stage corrections are made to the state predictions via noisy measurements.  A summary of the EKF (Ristic et. al \cite{ristic2004beyond}) equations is given below.  The target state $\mathbf{x}_k$ and measurement $\mathbf{z}_k$ equations propagate according to:
\begin{eqnarray}
\mathbf{x}_k &=& \mathbf{f}_{k-1}\left(\mathbf{x}_{k-1}\right) + \mathbf{v}_{k-1},\\
\mathbf{z}_k &=& \mathbf{h}_{k}\left(\mathbf{x}_{k}\right) + \mathbf{w}_{k},
\end{eqnarray}
where $\mathbf{v}_{k-1}$ and $\mathbf{w}_{k}$ are random sequences and are mutually independent with zero-mean, white Gaussian with covariances $\mathbf{Q}_{k-1}$ and $\mathbf{R}_{k}$ respectively.  The EKF is based on the assumption that local linearisation of the above equations may be a sufficient description of nonlinearity.  The mean and covariance of the underlying Gaussian density are computed recursively in a two stage process (Ristic et. al \cite{ristic2004beyond}):\\

Stage 1: Prediction
\begin{eqnarray}
\mathbf{\hat{x}}_{k|k-1} &=& \mathbf{f}\left(\hat{x}_{k-1|k-1}\right),\\ 
\mathbf{P}_{k|k-1} &=& \mathbf{Q}_{k-1} + \hat{\mathbf{F}}_{k-1}\hat{\mathbf{P}}_{k-1|k-1}\hat{\mathbf{F}}_{k-1}^\top.
\end{eqnarray}  

Stage 2: Update/Correction
\begin{eqnarray}
\mathbf{\hat{x}}_{k|k} &=& \mathbf{\hat{x}}_{k|k-1} + \mathbf{K}_k \mathbf{\nu}_k, \\
\mathbf{P}_{k|k} &=& \mathbf{P}_{k|k-1} - \mathbf{K}_k\mathbf{S}_k\mathbf{K}_k^\top,
\end{eqnarray}

where $\mathbf{\nu}_k = \mathbf{z}_k-\mathbf{h}_k\left(\mathbf{\hat{x}}_{k|k-1}\right)$, $\mathbf{S} = \mathbf{\hat{H}}_k\mathbf{P}_{k|k-1}\mathbf{\hat{H}}_k^\top + \mathbf{R}_k$ and $\mathbf{K}_k = \mathbf{P}_{k|k-1}\mathbf{\hat{H}}_k^\top\mathbf{S}_k^{-1}$.  The parameter $\mathbf{K}_k$ is commonly known as the Kalman gain and $\mathbf{S_k}$ is referred to as the innovation covariance.   The innovation $\mathbf{\nu}_k$ is the error between the predicted measurement and the actual measurement of the system.  The matrices $\mathbf{\hat{F}}_{k-1}$ and $\mathbf{\hat{H}}_k$ are the local linearisation of the nonlinear functions $\mathbf{f}_{k-1}$ and $\mathbf{h}_{k}$.  The two matrices are defined as Jacobians evaluated at $\mathbf{\hat{x}}_{k-1|k-1}$ and $\mathbf{\hat{x}}_{k|k-1}$ respectively (Ristic et. al \cite{ristic2004beyond}).  The non-Gaussianity of the true posterior density is more evident, for example becomes bimodal or heavily skewed, when the model is highly nonlinear.  In this event the performance of the EKF will significantly degrade.

\subsection{The Unscented Kalman Filter (UKF)} \label{subsec:chap4_UKF_equations}
The UKF addresses the issue of non-Gaussianity.  The UKF is a part of a family of nonlinear filters, referred to as linear regression Kalman Filters, that are based on statistical linearisation rather than analytical linearisation.  The key concept behind these filters is to perform nonlinear filtering using a Gaussian representation of the posterior $p\left(\mathbf{x}_k\vert\mathbf{Z}_k\right)$ through a set of deterministically chosen sample points.  The true mean and covariance of the Gaussian density are completely captured by these sample points up to the second order of nonlinearity, with errors introduced in the third and higher order when propagated through a nonlinear transform.  The EKF on the other hand is only of first order with errors introduced in the second and higher orders.  The filters belonging to this family differ only by the method used to select the sample points i.e. their number, weights and values in the filtering equations are identical and are given below.  The UKF uses an unscented transform for the selection of points in an EKF framework (Ristic et. al \cite{ristic2004beyond}).\\

We assume that at time $k-1$ the posterior is Gaussian: $p\left(\mathbf{x}_{k-1}\vert\mathbf{Z}_{k-1}\right) \approx \mathcal{N}\left(\mathbf{x}_{k-1};\mathbf{\hat{x}}_{k-1\vert{k-1},\mathbf{P}_{k-1\vert k-1}}\right)$.  The very first step is representing this density via a set of $N$ sample points $\mathcal{X}_{k-1}^i$ and their weights $W_{k-1}^i,\,i = 0,\hdots,N-1$.  The UKF uses the unscented transform \cite{ristic2004beyond} to select the sample points $\mathcal{X}_{k-1}^i$ and weights $W_{k-1}^i$.  The prediction step is as follows:
\begin{eqnarray}
\mathbf{\hat{x}}_{k\vert k-1} &=& \sum_{i = 1}^{N-1}W_{k-1}^i\mathbf{f}_{k-1}\left(\mathcal{X}_{k-1}^i\right),\\
\mathbf{P}_{k\vert k-1} &=& \mathbf{Q}_{k-1} + \sum_{i=0}^{N-1}W_{k-1}^i\left[\mathbf{f}_{k-1}\left(\mathcal{X}_{k-1}^i\right)-\mathbf{\hat{x}}_{k\vert k-1}\right]\left[\mathbf{f}_{k-1}\left(\mathcal{X}_{k-1}^i\right) - \mathbf{\hat{x}}_{k\vert k -1}\right]^\top.
\end{eqnarray} 
A set of $N$ sample points:
\begin{equation}
\mathcal{X}_{k\vert k-1}^i = \mathbf{f}_{k-1}\left(\mathcal{X}_{k-1}^i\right),
\end{equation}
are used to represent the predicted density: $p\left(\mathbf{x}_k\vert \mathbf{Z}_{k-1}\right) \approx \mathcal{N}\left(\mathbf{x}_k;\mathbf{\hat{x}}_{k\vert k-1},\mathbf{P}_{k \vert k-1}\right)$ and the predicted measurement becomes:
\begin{equation}
\mathbf{\hat{z}}_{k \vert k-1} = \sum_{i = 0}^{N-1} W_{k-1}^i\mathbf{h}\left(\mathcal{X}_{k-1}^i\right).
\end{equation}
The update step is defined as:
\begin{eqnarray}
\mathbf{\hat{x}}_{k|k-1} &=& \mathbf{\hat{x}}_{k \vert k-1} + \mathbf{K}_k\nu_k,\\
\mathbf{P}_{k \vert k} &=& \mathbf{P}_{k \vert k-1} - \mathbf{K}_k\mathbf{S}_k\mathbf{K}_k^\intercal,  
\end{eqnarray}

where:
\begin{eqnarray}
\mathbf{K}_k &=& \mathbf{P}_{xz}\mathbf{S}_k^{-1},\\
\nonumber \\ 
\mathbf{S}_k &=& \mathbf{R}_k + \mathbf{P}_{zz},\\
\nonumber \\
\mathbf{P}_{xz} &=& \sum_{i=0}^{N-1}W_{k-1}^i\left(\mathcal{X}_{k\vert k-1}^i- \mathbf{\hat{x}}_{k\vert k-1}\right)\left(\mathbf{h}_k(\mathcal{X}_{k\vert k-1}^i) - \mathbf{\hat{z}}_{k\vert k-1}\right)^{\intercal},\\
\nonumber \\
\mathbf{P}_{zz} &=& \sum_{i=0}^{N-1}W_{k-1}^i\left(\mathbf{h}_k(\mathcal{X}_{k\vert k-1}^i) - \mathbf{\hat{z}}_{k\vert k-1}\right)\left(\mathbf{h}_k(\mathcal{X}_{k\vert k-1}^i) - \mathbf{\hat{z}}_{k\vert k-1}\right)^{\intercal}
\end{eqnarray}

As can be seen from the above filter equations there is no explicit calculation of Jacobians.  Consequently these filters can be utilised even when the nonlinear functions $\mathbf{f}$ and $\mathbf{h}$ have discontinuities.

\subsection{Interacting Multiple Model (IMM)}\label{subsection:chap4_IMM_equations}
The IMM belongs to a class of filters called the Gaussian Sum Filters.  The main concept here is the approximation of the required posterior density $p\left(\mathbf{x}_k\vert \mathbf{Z}_k\right)$ by a Gaussian mixture (Ristic et. al \cite{ristic2004beyond}):
\begin{equation}
p\left(\mathbf{x}_k\vert \mathbf{Z}_k\right)\approx p_A\left(\mathbf{x}_k\vert \mathbf{Z}_k\right) = \sum_{i=1}^{q_k}w_k^i\,\,\mathcal{N}\left(\mathbf{x}_k^i;\mathbf{\hat{x}}_{k\vert k}^i,\mathbf{P}_{k\vert k}^i\right),
\end{equation}
where $w_k^i$ are weights that are normalised, $\sum_{i=1}^{q_k} w_k^i = 1$.  Gaussian sum filters such as the IMM are ideal when the posterior density is multimodal because for multimodal densities there is a performance degradation in both the EKF and UKF.  At time $k$ the state estimate is calculated for each possible current model using $r$ filters, with each filter using a different combination of the previous model-conditioned estimates called mixed initial condition.  The algorithm as outlined in Bar-Shalom et. al \cite{bar2004estimation} is:
\begin{itemize}
\item[Step 1:] Calculation of the mixing probabilities.  The probability that mode $M_i$ was in effect at time $k-1$ given that $M_j$ is in effect at $k$ conditioned on $Z^{k-1}$ is given by:
\begin{align}
\mu_{i,j}\left(k-1 \vert k-1 \right) &\triangleq P\left\lbrace M_i\left(k-1\right)\vert M_j\left(k\right),Z^{k-1}\right\rbrace, \\
&= \frac{1}{\bar{c}_j}P\left\lbrace M_j\left(k\right)\vert M_i\left(k-1\right),Z^{k-1}\right\rbrace P\left\lbrace M_i \left(k-1\right)\vert Z^{k-1}\right\rbrace.
\end{align}
The above can be written as:
\begin{equation}
\mu_{i,j}\left(k-1 \vert k-1 \right) = \frac{1}{\bar{c}_j}p_{ij}\,\mu_i\left(k-1\right),\quad\quad i,j = 1,\hdots,r,
\end{equation}
where the normalising constants are:
\begin{equation}
\bar{c}_j=\sum_{i=1}^{r}p_{ij}\,\mu_i\left(k-1\right),\quad\quad i,j = 1,\hdots,r.
\end{equation}
\item[Step 2:] Mixing.  The mixed initial condition for the filter matched to $M_j\left(k\right)$ is calculated starting with $\hat{x}^i\left(k-1\vert k-1\right)$:
\begin{equation}
\hat{x}^{0j}\left(k-1\vert k-1\right) = \sum_{i = 1}^{r}\hat{x}^i\left(k-1\vert k-1\right) \mu_{i\vert j} \left(k-1\vert k-1 \right),\quad\quad i,j = 1,\hdots,r,
\end{equation}
and the corresponding covariance is given by:
\begin{equation}
\begin{split}
P^{0j}\left(k-1\vert k-1\right) = &\sum^{r}_{i=1} \mu_{i\vert j}\left(k-1\vert k-1\right)\bigg\{ P^i\left(k-1\vert k-1\right)\\
&+ \left[\hat{x}^i\left(k-1\vert k-1\right)-\hat{x}^{0j}\left(k-1\vert k-1\right)\right]\\
&\cdot \left[\hat{x}^i\left(k-1 \vert k-1\right) - \hat{x}^{0}\left(k-1\vert k-1\right)\right]^\intercal\bigg\},\quad\quad i,j = 1,\hdots,r.
\end{split}
\end{equation}
\item[Step 3:] Mode-Matched Filtering.  The estimates of the states and covariances calculated in step 2 above are used as inputs to the filter matched to $M_j\left(k\right)$ which uses $z\left(k\right)$ to determine $\hat{x}^j\left(k\vert k\right)$ and $P^j\left(k\vert k\right)$.  The likelihood functions associated to the $r$ filters:
\begin{equation}
\Lambda_j\left(k\right)=p\left[z\left(k\right)\vert M_j\left(k\right),Z^{k-1}\right],
\end{equation}
are calculated using the mixed initial condition and covariance from step 2:
\begin{equation}
\Lambda_j\left(k\right) = p\left[z\left(k\right)\vert M_j\left(k\right),\hat{x}^{0j}\left(k-1\vert k-1\right),P^{0j}\left(k-1\vert k-1\right)\right],\quad\quad j = 1,\hdots,r,
\end{equation}
that is:
\begin{equation}
\Lambda_j\left(k\right) = \mathcal{N}\left[z\left(k\right);\hat{z}^j\left[k\vert k-1;\hat{x}^{0j}\left(k-1\vert k-1\right)\right],S^j\left[k;P^{0j}\left(k-1\vert k-1\right)\right]\right],\quad\quad j = 1,\hdots,r.
\end{equation}
\item[Step 4:] Mode Probability Update.  The mode probabilities are then updated via:
\begin{equation}
\mu_{j}\left(k\right)=\frac{1}{c}\Lambda_{j}\left(k\right)\bar{c}_j,
\end{equation}
where $c$ is the normalisation constant and is given by $c = \sum_{j=1}^r\Lambda_{j}\left(k\right)\bar{c}_{j}$
\item[Step 5:] Estimations and Covariance Combination.  The output is obtained by combining the model-conditioned estimates and covariances:
\begin{eqnarray}
\hat{x}\left(k\vert k\right) &=& \sum_{j=1}^r\hat{x}^j\left(k\vert k\right)\mu_j\left(k\right),\\
\nonumber \\
\hat{P}\left(k\vert k\right) &=& \sum_{j=1}^r \mu_j\left(k\right)\left\lbrace P^j\left(k\vert k\right)+\left[\hat{x}^j\left(k\vert k\right) - \hat{x}\left(k\vert k\right)\right]\left[\hat{x}^j\left(k\vert  k\right) - \hat{x}\left(k\vert k\right)\right]^\intercal\right\rbrace.
\end{eqnarray} 
\end{itemize} 

This section outlined the details of the the methods chosen for further investigation, the EKF, the UKF and the IMM filters.  These techniques are applied to the 2D robot model in the next section. 

\section{Problem Formulation} \label{section:chap4_ProbForm}
\begin{figure}[H]
\begin{center}
\includegraphics[scale=0.5]{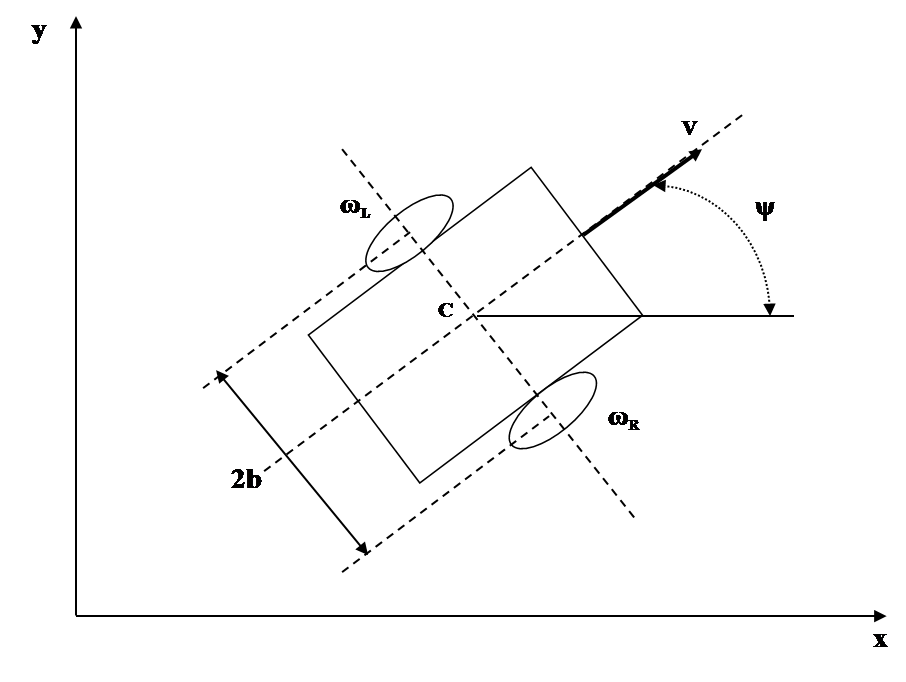} 
\caption{Robot Schematic}
\label{fig:chap3_Robot EOM}
\end{center}
\end{figure}

To test the different filtering techniques in an FDI context the robot model of \cite{RKPaper1} is used for development and testing purposes and forms the plant that is to be controlled as well as the process model of the NMPC controller (see figure \ref{fig:chap3_Robot EOM}).

The equations for this robot model are:
\begin{eqnarray}
\dot{x} &=&  V\cos\psi, \\
\dot{y} &=&  V\sin\psi, \\
\dot{\psi} &=& \frac{R(\omega_R - \omega_L)}{2b},
\end{eqnarray}
where $x$ is the $x$-coordinate of the point $C$, $y$ is the $y$-coordinate of the point $C$, $\psi$ is the heading angle, $\omega_R$ is the right wheel angular velocity, $\omega_L$ is the left wheel angular velocity and $V$ is the speed given by $V = \frac{R(\omega_R + \omega_L)}{2}$.\\

The NMPC controller is based on the design given in \cite{RKPaper1} which solves the following optimal control problem using pseudospectral discretisation:

\begin{equation}
\min \, \frac{\left(t_f-t_0\right)}{2}\sum_{j=0}^N\,\left(\big\Vert\mathbf{x}-\mathbf{x}_{\text{ref}}\big\Vert_{Q_x}^2+\big\Vert V-V_{\text{ref}}\big\Vert_{Q_V}^2 + \big\Vert \dot{\psi}-\dot{\psi}_{\text{ref}}\big\Vert_{Q_\psi}^2\right)\,w_j.
\end{equation}

subject to
\begin{eqnarray}
\left(\frac{t_f-t_0}{2}\right)\mathbf{D}_{j,k}\mathbf{x}_j - \mathbf{\dot{x}}_j &=& 0, \\
\mathbf{x}(j_0) - \mathbf{x}_{\text{dem}}(j_0) &=& 0,\\
\mathbf{x}_{lb}  \leq    \mathbf{x}  \leq  \: \mathbf{x}_{ub},\\
\mathbf{u}_{lb}  \leq    \mathbf{u}  \leq  \: \mathbf{u}_{ub},\\
\Delta\mathbf{\delta}_{e_{\text{lb}}}  \leq    \Delta\mathbf{\delta}_e  \leq  \: \Delta\mathbf{\delta}_{e_{\text{ub}}}, \label{eq:chap5_3DOF_cons}
\end{eqnarray}

Full details of the design are given in \cite{RKPaper1}.  Fifty coincidence points ($N = 50$) are used along with a prediction window length of 5 secs ($H_p = 5$).  The weights $Q_x,\, Q_V,$ and $Q_\psi$  are diagonal matrices with the diagonal values set to 10, 1 and 1 respectively and found through trial and error.\\

The fault, to be simulated and tested for, is a punctured wheel.  If a wheel is punctured the radius of the wheel will decrease and so the filters are set up to estimate the radius of the wheel.  Four different filters have been designed, the EKF, the UKF, the EKF IMM and the UKF IMM.\\

The robot parameters used for all simulations are: Right wheel radius, $R_R=2\,m$, Left wheel radius, $R_L=2\,m$, Distance between wheels, $b=1\,m$, Speed demand is $10\,m/s$ and the input constraints on $\omega_R$ and $\omega_L$ are $\pm 1000\,\deg/\sec$.  The filters are updated at 100Hz while the controller is updated at 10Hz.  All work was developed using MATLAB with SNOPT as the nonlinear programme (NLP) solver.  The following subsections detail the design of each of the filters.

\subsection{EKF Fault Detection Filter}\label{subsec:chap4_EKF_equations_sims}
The state vector for the EKF consists of the following states:
\begin{equation} \label{eqn:chap4_EKF_stateVec}
\mathbf{x} = \left[x,\,\,\, y,\,\, \psi,\,\,\, R_R,\,\,\, R_L\right]^\intercal,
\end{equation}
where $x$, $y$ and $\psi$ are the robot states and $R_R$ and $R_L$ are the right wheel and left wheel radii respectively.  The measurements are assumed to be of the speed, $V$, of the robot:

\begin{equation} \label{eqn:chap4_EKF_zVec}
z = V(k) + w(k),
\end{equation}
where $\mathbf{w}(k)$ is additive white noise.  The initial state vector and initial covariance matrix are:
\begin{equation} \label{eqn:chap4_EKF_xPInit}
\mathbf{x}(0) = \left[x_0,\,\,\, y_0,\,\, \psi_0,\,\,\, 2,\,\,\, 2\right]^\intercal,\quad\quad
P(0) = \begin{bmatrix}
    (0.5)^2 & 0 & 0 & 0 & 0 \\
    0 & (0.5)^2 & 0 & 0 & 0 \\
    0 & 0 & (1\pi/180)^2 & 0 & 0 \\
    0 & 0 & 0 & (0.5)^2 & 0 \\
    0 & 0 & 0 & 0 & (0.5)^2 
\end{bmatrix}.
\end{equation}
The $Q$ and $R$ noise matrices were chosen to be:
\begin{equation} \label{eqn:chap4_EKF_QR}
Q = \begin{bmatrix}
    (5\,\Delta t)^2 & 0 & 0 & 0 & 0 \\
    0 & (5\,\Delta t)^2; & 0 & 0 & 0 \\
    0 & 0 & (0.1\,\Delta t)^2 & 0 & 0 \\
    0 & 0 & 0 & (2\,\Delta t)^2 & 0 \\
    0 & 0 & 0 & 0 & (2\,\Delta t)^2 
\end{bmatrix},\quad\quad
R = (0.5)^2,
\end{equation}
where $\Delta t$ is the update rate of the filter.  For the prediction cycle an Euler integration scheme is used to predict the states of the EKF forward. The predicted measurement $\hat{z}$ is given by:
\begin{equation}
\hat{z}(k) = V(k|k-1),
\end{equation}
Given the above information the Kalman filter equations given in section \ref{subsection:chap4_EKF_equations} are applied to estimate the radius of each wheel in the experiments conducted in section \ref{section:chap4_NRandA}. 

\subsection{UKF Fault Detection Filter}
The general structure of the EKF and UKF are very similar in that they both have a prediction and update cycle and produce a single state vector and a corresponding covariance matrix.  For the robot model the state vector is the same as the one given in equation \eqref{eqn:chap4_EKF_stateVec}.  The initial state vector, initial covariance matrix, the process noise matrix $Q$ and the noise covariance matrix $R$ all remain the same as those given in subsection \ref{subsec:chap4_EKF_equations_sims}.  The UKF algorithm given in subsection \ref{subsec:chap4_UKF_equations} is applied to the robot model with $\kappa = 0.001$ \cite{JulierAndUhlmann}.

\subsection{Interacting Multiple Model Fault Detection Filter}
The interacting multiple model method, as the name suggests, is made up of multiple models where each model tests a different hypothesis.  Four different models (the terms mode and model are used interchangeably and have the same meaning in the context of IMMs) have been designed where:

\begin{itemize}[leftmargin=0.65in]
\item[Model 1:] No Fault case.
\item[Model 2:] $50\%$ right wheel deflation, left wheel no fault.
\item[Model 3:] Right wheel no fault, $50\%$ left wheel deflation.
\item[Model 4:] $50\%$ right wheel deflation, $50\%$ left wheel deflation.
\end{itemize}  

During Step 3 of the IMM algorithm given in section \ref{subsection:chap4_IMM_equations} a filter such as the EKF is used to update the states and covariances, and both an EKF based IMM filter and a UKF based IMM filter have been designed.  The initial covariance matrix for each filter and each mode are the same as equation \eqref{eqn:chap4_EKF_xPInit}.  The $Q$ and $R$ matrices are those given in equation \eqref{eqn:chap4_EKF_QR} and the initial state vectors for each filter and mode are:
\begin{eqnarray}
\mathbf{x}_1(0) &=& \left[x_0,\,\,\,y_0\,\,\,\psi_0\,\,\,2\,\,\,2\right]^\intercal,\\
\mathbf{x}_2(0) &=& \left[x_0,\,\,\,y_0\,\,\,\psi_0\,\,\,1\,\,\,2\right]^\intercal,\\
\mathbf{x}_3(0) &=& \left[x_0,\,\,\,y_0\,\,\,\psi_0\,\,\,2\,\,\,1\right]^\intercal,\\
\mathbf{x}_4(0) &=& \left[x_0,\,\,\,y_0\,\,\,\psi_0\,\,\,1\,\,\,1\right]^\intercal.
\end{eqnarray}
The mixing probabilities or mode probabilities are initially set to:
\begin{equation}
\mu = \left[1/4,1/4,1/4,1/4\right]^\intercal,
\end{equation}
and the mode transition probabilities matrix $p$ is set to:
\begin{equation}
p = \begin{bmatrix}
0.97 & 0.01 & 0.01 & 0.01\\
0.01 & 0.97 & 0.01 & 0.01\\
0.01 & 0.01 & 0.97 & 0.01\\
0.01 & 0.01 & 0.01 & 0.97
\end{bmatrix}.
\end{equation} 

\section{Numerical Results and Analysis}\label{section:chap4_NRandA}
The robot is required to follow a circular path for all experiments.  To simulate the measurement additive white noise is added to the speed of the robot which is calculated as a part of the truth simulations of the robot movement.  To test the filters four different test cases were set up and each test case was run twice.  During the first run the FDI feedback loop is not closed and the filters are used for estimation only.  The FDI loop is closed during the second run to investigate the behaviour of the full active FTC controller.  The test cases are as follows: \\

\begin{enumerate}[leftmargin=2cm, label=\bfseries Scenario \arabic*:]
\item No Fault.  The objective is to investigate how well the filters estimate the radii of the tyres in a no fault situation.
\item Left wheel $50\%$ puncture.  In this case a puncture is simulated to occur 10 secs into the simulation.  The wheel is assumed to deflate to $50\%$ of its original value instantaneously.
\item Left and Right Wheel puncture.  In this test case a left wheel puncture is simulated 5 secs into the run and a right wheel puncture is simulated to occur 10 secs into the simulation.  Both punctures are assumed to cause an instantaneous reduction of the respective wheel radius to $50\%$ of the original wheel radius.
\item Left wheel linear puncture.  In this test case once again the left wheel is punctured 10 secs into the run however this time the puncture is assumed to follow a linear reduction in wheel radius according to $R_L = 2 - 0.1t$, where $R_L$ represents the left wheel radius reduced from its original value of 2m down to 0.1m at a rate of 0.1m/s and t is the current time.  The radius does not drop to 0m as this caused a complete system failure.\\   
\end{enumerate}
The results for each filter are presented in the next four subsections:

\subsection{Scenario 1}
Plots of the speed innovation were produced (but have been omitted due to space constraints) where the innovations were plotted along with the calculated $2\sigma$ uncertainty bounds.  The $2\sigma$ uncertainty bounds are a $95\%$ confidence interval and the solution (innovations in this case) must remain within the bounds $95\%$ of the time.  The results for the EKF and UKF filters showed that the speed innovations remained well within the uncertainty bounds throughout the duration of the run both with and without feedback.  The speed innovation plots for the IMM filters were also produced.  The results for both the UKF and EKF based IMM filters showed that the filters were able to very quickly detect the correct mode of operation.\\  

Plots of the wheel radius estimates were also produced.  All filters do a very good job of estimating the radius of the tyres with and without feedback.  The IMM filters initially have a higher error in the tyre estimate as the estimation is based on a mixture of all the models, however it took only one update for the IMM to reach the correct estimate.\\

The wheel speed estimates (or angular rates) were also analysed and it was found that the estimates for all four filters were very similar.  In the case of no feedback the estimate was the same as the actual speed; however in the case where feedback is provided the wheel speeds were quite noisy.  This is a result of calculations based on noisy measurements which is a consequence of feedback.\\

Plots of the robot trajectory for all four filters showed that the robot remained on the path with and without feedback which is to be expected in the no fault case.

\subsection{Scenario 2}
The speed innovations for all four filters were plotted and the results showed that all filters were able to detect the fault.  The fault occurred at 10 secs into the run, and the plots show that at 10 seconds there was a peak change in the innovation curves.  The corrections/innovations were seen to increase at the time the fault occurred and settle again near zero once the correct estimate was reached.\\

\begin{figure}[H]
\center
\includegraphics[scale=0.3]{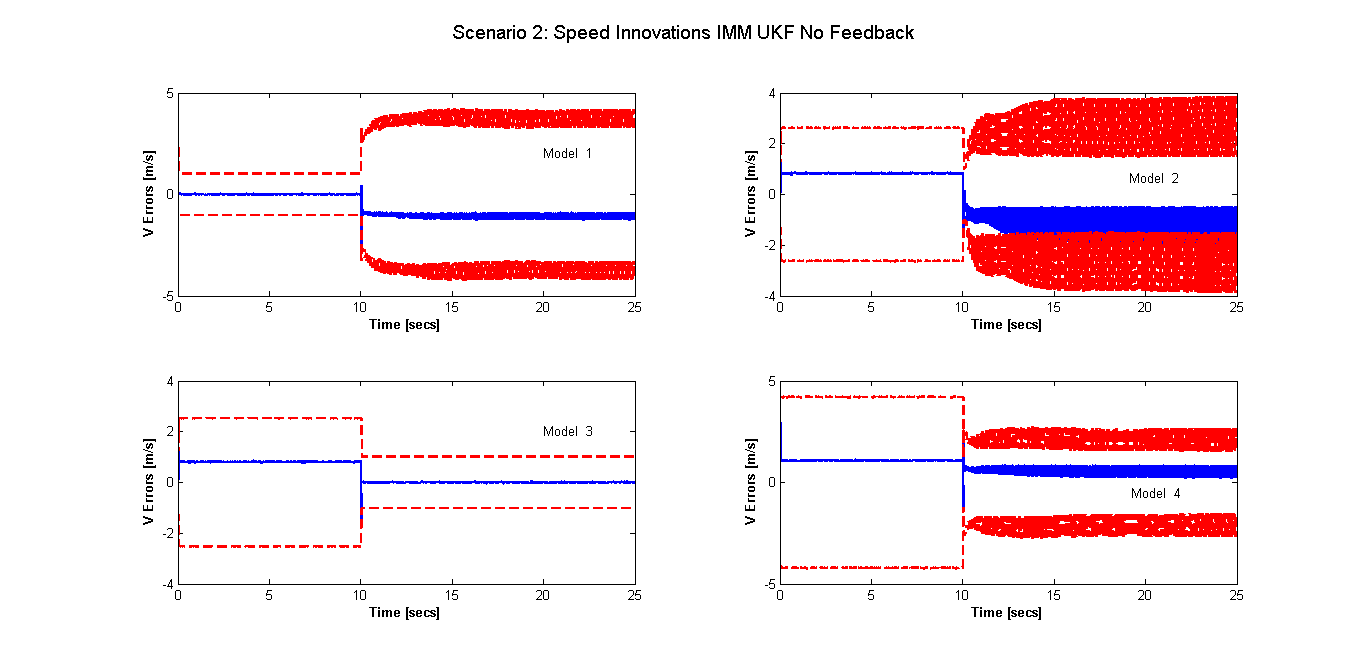}
\caption{Scenario 2 - IMM UKF Speed Innovations No Feedback, $2\sigma$ Uncertainty Bounds (red), Speed Innovations (blue)} 
\label{fig:chap4_scenario2_IMMUKF_velocityErrors_NFB}
\end{figure}

\begin{figure}[H]
\center
\includegraphics[scale=0.3]{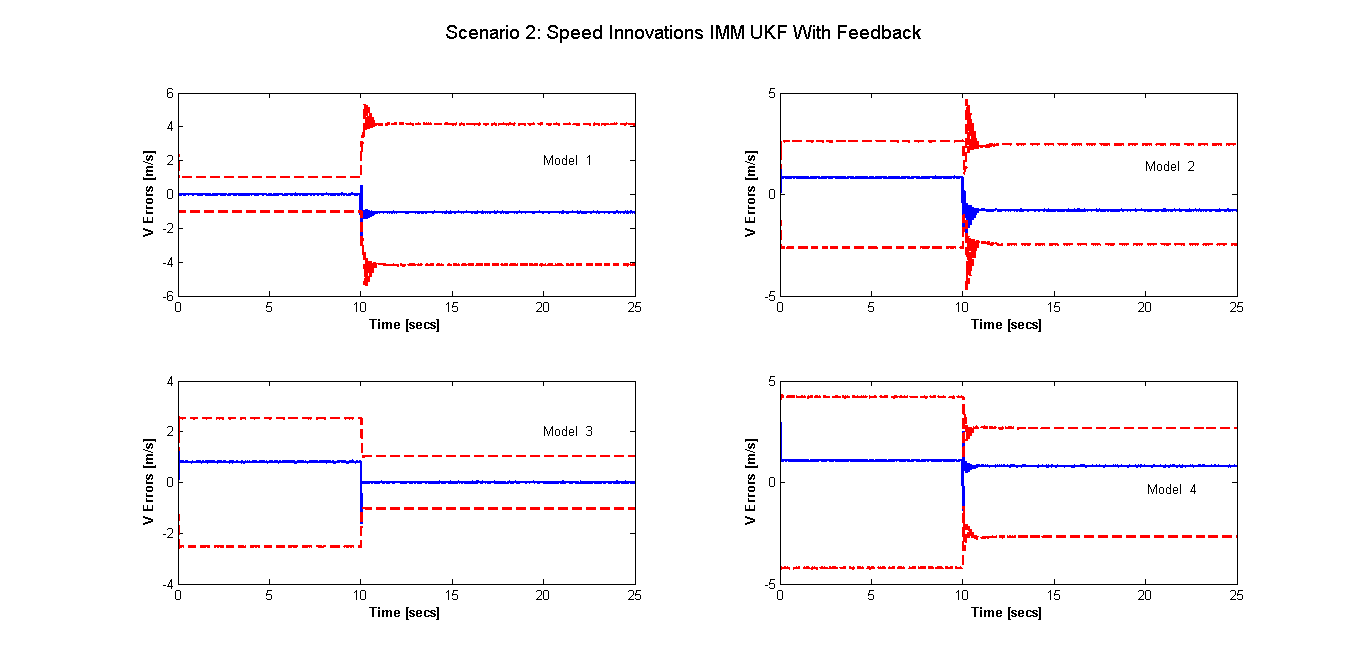}
\caption{Scenario 2 - IMM UKF Speed Innovations With Feedback, $2\sigma$ Uncertainty Bounds (red), Speed Innovations (blue)} 
\label{fig:chap4_scenario2_IMMUKF_velocityErrors_WFB}
\end{figure}

For the IMM filters model 3 is the correct match for scenario 2 and both of the IMM filters were seen to find the correct mode immediately as mode 3 is the most confident in its estimate.  Mode 3 and Mode 4 both hypothesise a failure in the left wheel of $50\%$ which is why after the occurrence of the fault the uncertainties do not increase.  However because mode 1 and mode 2 do not hypothesise a fault in the left wheel the uncertainties can be seen to increase once the fault has occurred.  The uncertainty in mode 3 was seen to decrease twice as much compared to mode 4 after the fault occurred.  This is because mode 4 hypothesises that both wheels are punctured whereas mode 3 predicts the puncture of only the left wheel.  Another point to note is that in the single filter cases feedback did not have much effect on the innovations.  However, in the case of the IMM filters the results show that with feedback the filter errors do not grow as rapidly between updates.  The errors are seen to grow very quickly when no feedback is present which is evident in figures \ref{fig:chap4_scenario2_IMMUKF_velocityErrors_NFB} and  \ref{fig:chap4_scenario2_IMMUKF_velocityErrors_WFB}.


The results for the wheel radius estimates showed that the UKF estimates were closer to the actual wheel radius compared to those produced by the EKF.  Turning on feedback results in the filters reaching a steady estimate faster when compared to the no feedback case.  The IMM filters produced slightly better estimates than the single UKF filter and there was very little improvement on the IMM estimates compared to the no feedback case.\\

The wheel speeds were also analysed and the results for all four filters presented the same trends.  Without feedback there is much more activity present compared to turning on the feedback.  Once the fault occurs the robot yaws to the side with the punctured wheel and demands a faster speed to compensate for the loss in radius.\\
 
\begin{figure}[H]
\center
\includegraphics[scale=0.3]{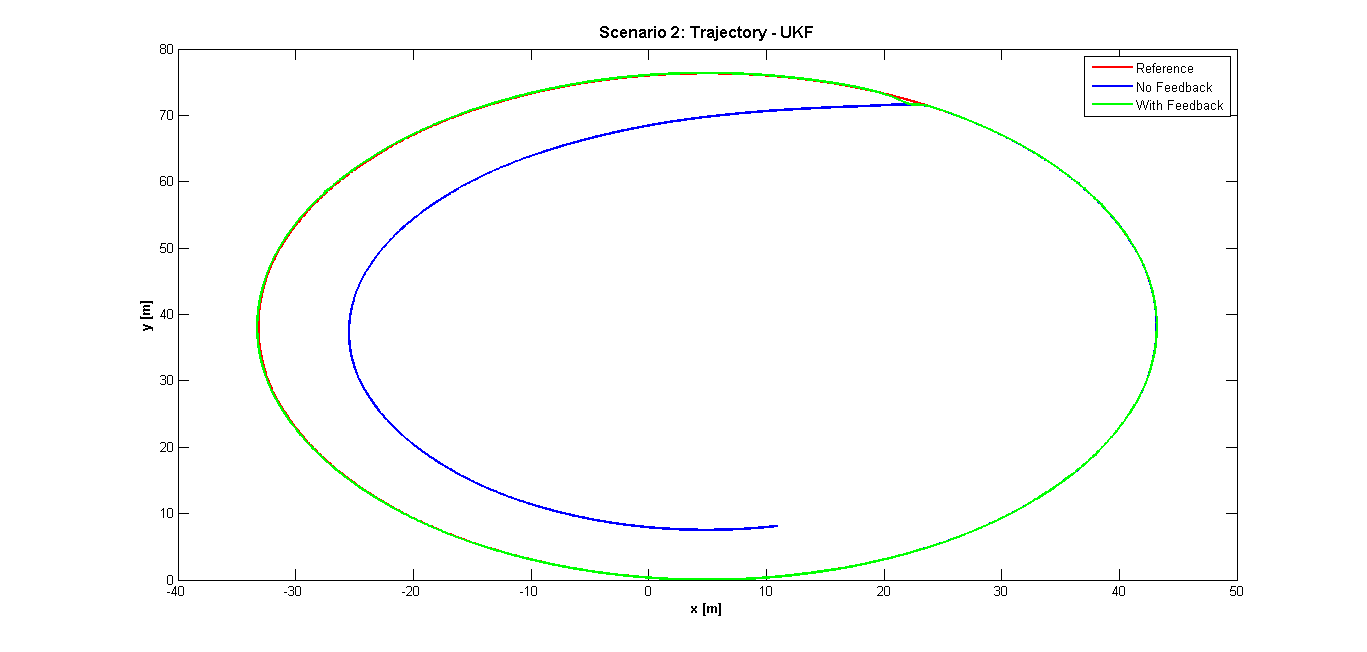} 
\caption{Scenario 2 - UKF Trajectory} 
\label{fig:chap4_scenario2_traj_UKF}
\end{figure}

The robot trajectories were simulated and all results showed that the robot was only able to remain on the path if feedback from the filter was provided to reconfigure the controller.  Although all filter estimates without feedback were excellent, without reconfiguration of the controller the robot could not be made to follow the desired path (see figure \ref{fig:chap4_scenario2_traj_UKF}).

\subsection{Scenario 3}
\begin{figure}[H]
\center
\includegraphics[scale=0.3]{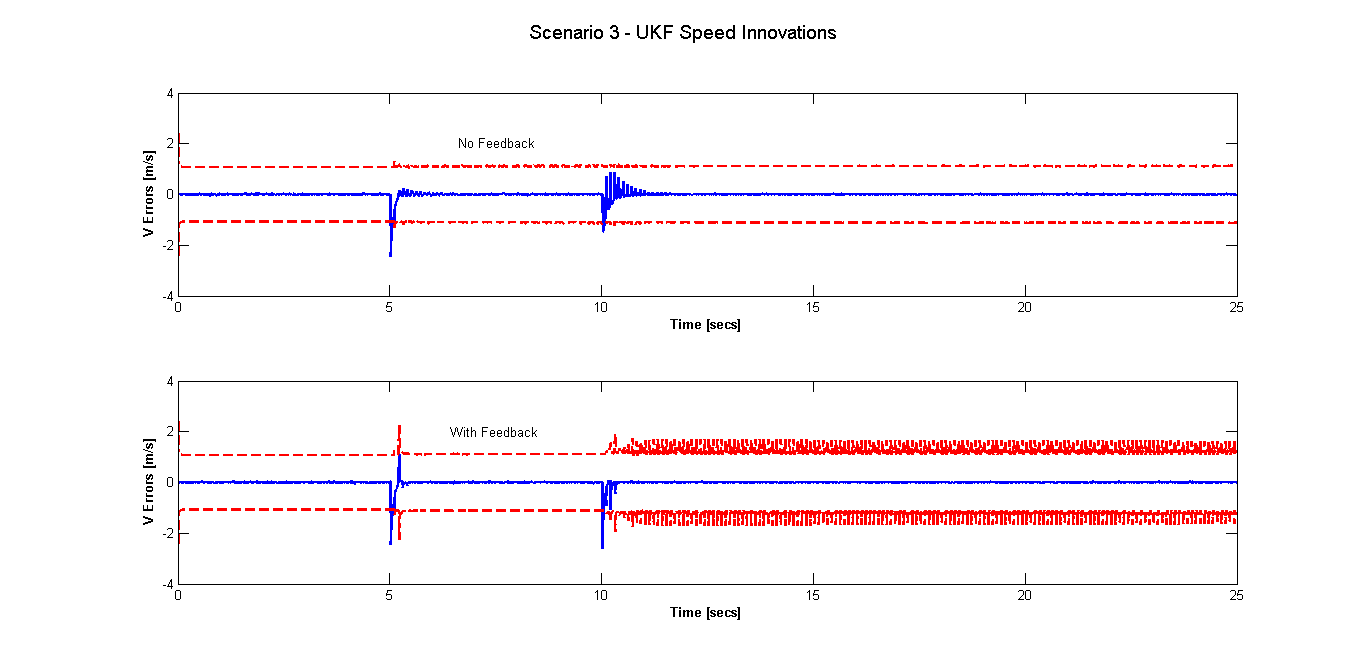}
\caption{Scenario 3 - UKF Speed Innovations, $2\sigma$ Uncertainty Bounds (red), Speed Innovations (blue)} 
\label{fig:chap4_scenario3_UKF_velocityErrors}
\end{figure}

\begin{figure}[H]
\center
\includegraphics[scale=0.3]{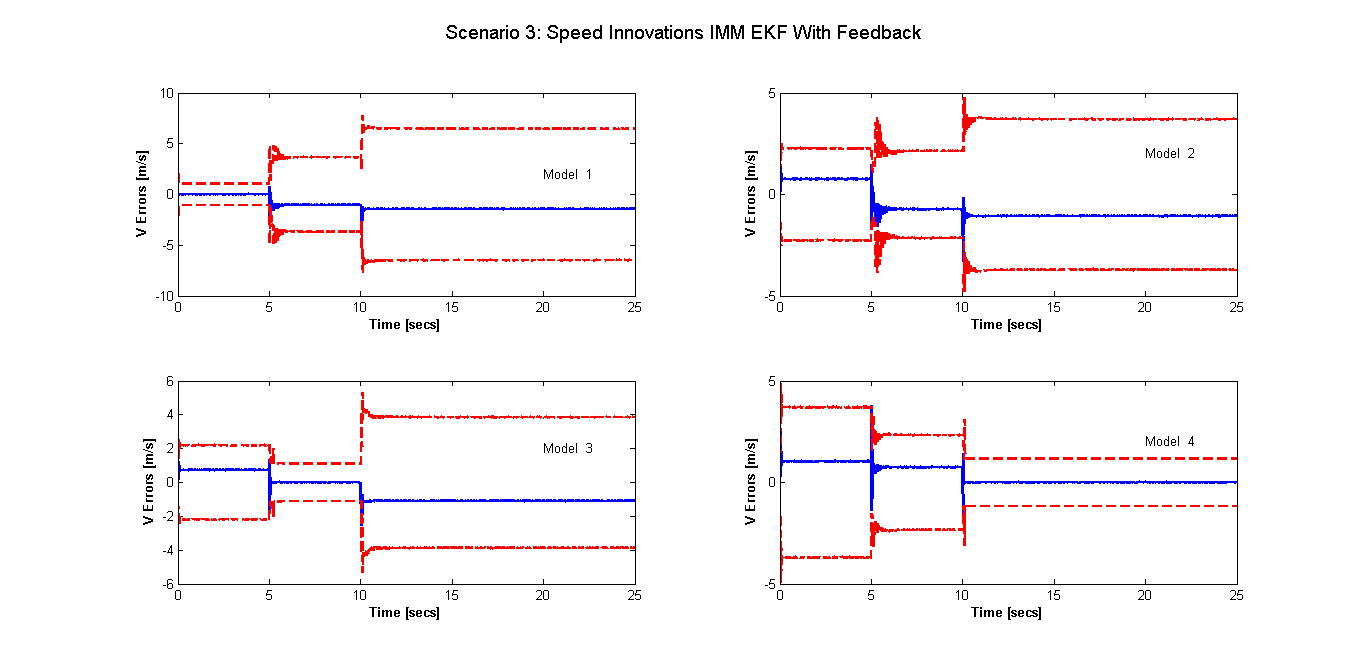}
\caption{Scenario 3 - IMM EKF Speed Innovations With Feedback, $2\sigma$ Uncertainty Bounds (red), Speed Innovations (blue)} 
\label{fig:chap4_scenario3_IMMEKF_velocityErrors_WFB}
\end{figure}

The plots of the speed innovations for all filters clearly indicated, from the sudden changes in innovations, the detection of both faults, left wheel at 5 secs and right wheel at 10 secs.  The EKF innovations were found to be consistent, however the innovations produced by the UKF show that, with feedback, the innovation uncertainty begins to grow rapidly between updates whereas without feedback the uncertainty remains constant (see figure \ref{fig:chap4_scenario3_UKF_velocityErrors}).  The IMM filters show that after 5 secs model 3 is the best model.  However, once the second fault occurs the filters do an excellent job of recognising that mode 4 is the correct match and uncertainties in mode 4 are seen to decrease (see figure \ref{fig:chap4_scenario3_IMMEKF_velocityErrors_WFB}).  It was again observed in the no feedback case that the uncertainties on the IMM filters grow rapidly between updates and many more corrections are required.\\   

\begin{figure}[H]
\center
\includegraphics[scale=0.3]{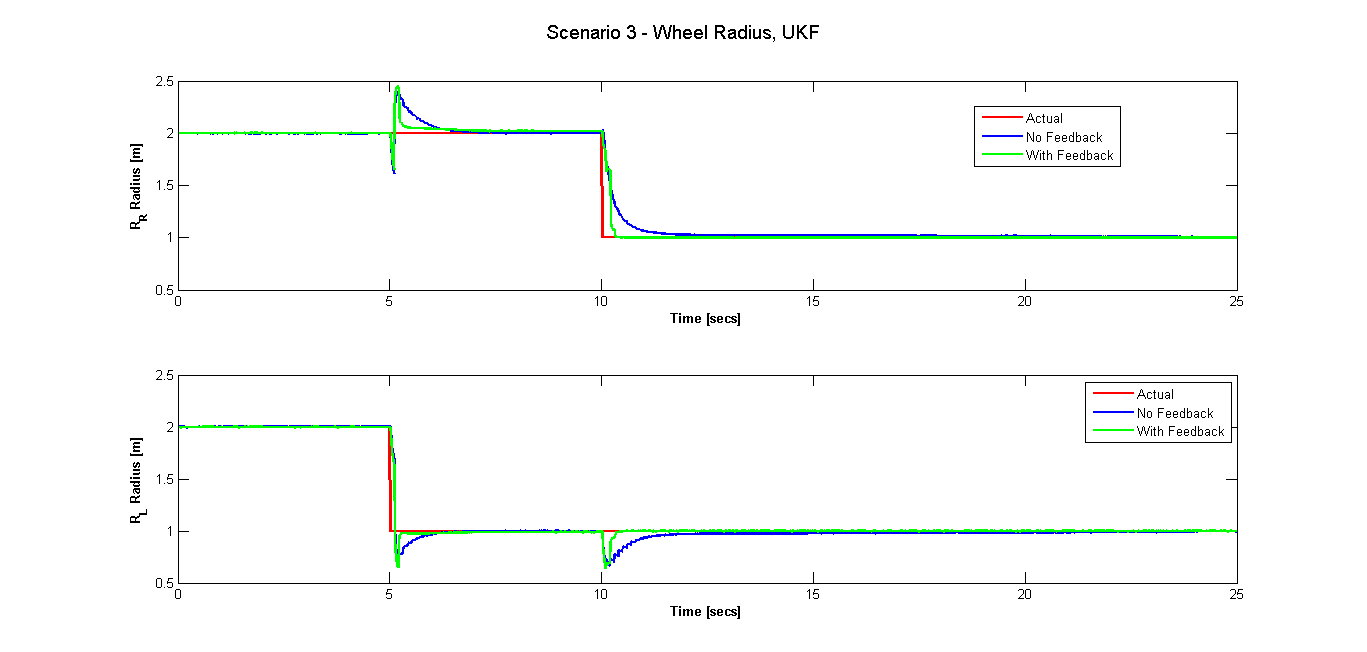}
\caption{Scenario 3 - UKF Radius Estimates, Right wheel (top), Left wheel (Bottom)} 
\label{fig:chap4_scenario3_UKF_wheelRadius}
\end{figure}

The wheel radius estimates produced by the filters showed that the IMM filters produce the best estimates of the radii.  The UKF performs slightly better than the EKF, and turning feedback on results in a faster settling time (see figure \ref{fig:chap4_scenario3_UKF_wheelRadius}).\\

Plots of the angular rates achieved by the robot via the EKF showed that once a fault occurred the punctured wheel is required to spin faster in order to compensate for the loss in radius.  The angular rates produced as a result of the UKF, with feedback, resulted in operation at the angular rate constraints.  Both IMM filters displayed similar behaviour in that once a wheel was punctured it was required to rotate faster to compensate for the loss in radius.\\

\begin{figure}[H]
\center
\includegraphics[scale=0.3]{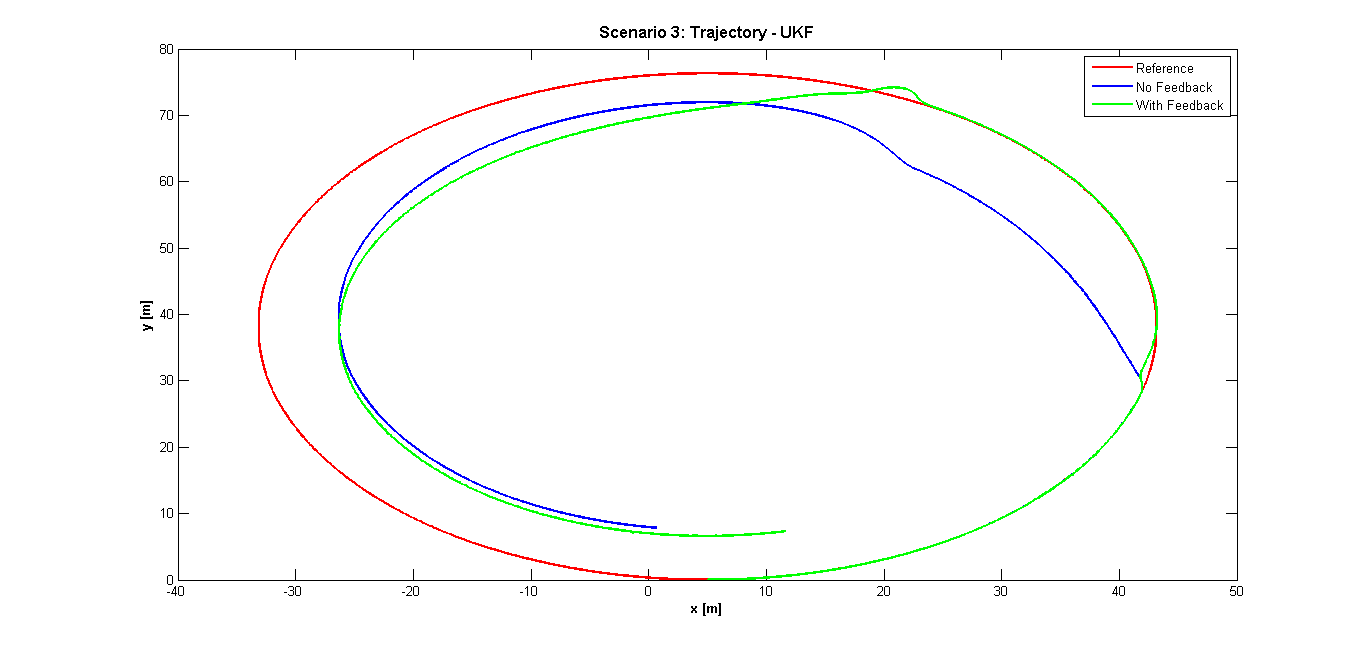} 
\caption{Scenario 3 - UKF Trajectory} 
\label{fig:chap4_scenario3_traj_UKF}
\end{figure}

The trajectories produced by each filter for scenario 3 showed that without feedback it is impossible to maintain the robot on the path.  Reconfiguring the controller on the other hand with estimates from the filters allowed the robot to easily follow the reference path.  An anomaly occurred with the UKF filter where even turning the feedback on did not result in the robot following the path after the occurrence of the second fault.  This could possibly be the result of poor tuning of the filter (see figure \ref{fig:chap4_scenario3_traj_UKF}).  \\

\subsection{Scenario 4}
\begin{figure}[H]
\center
\includegraphics[scale=0.3]{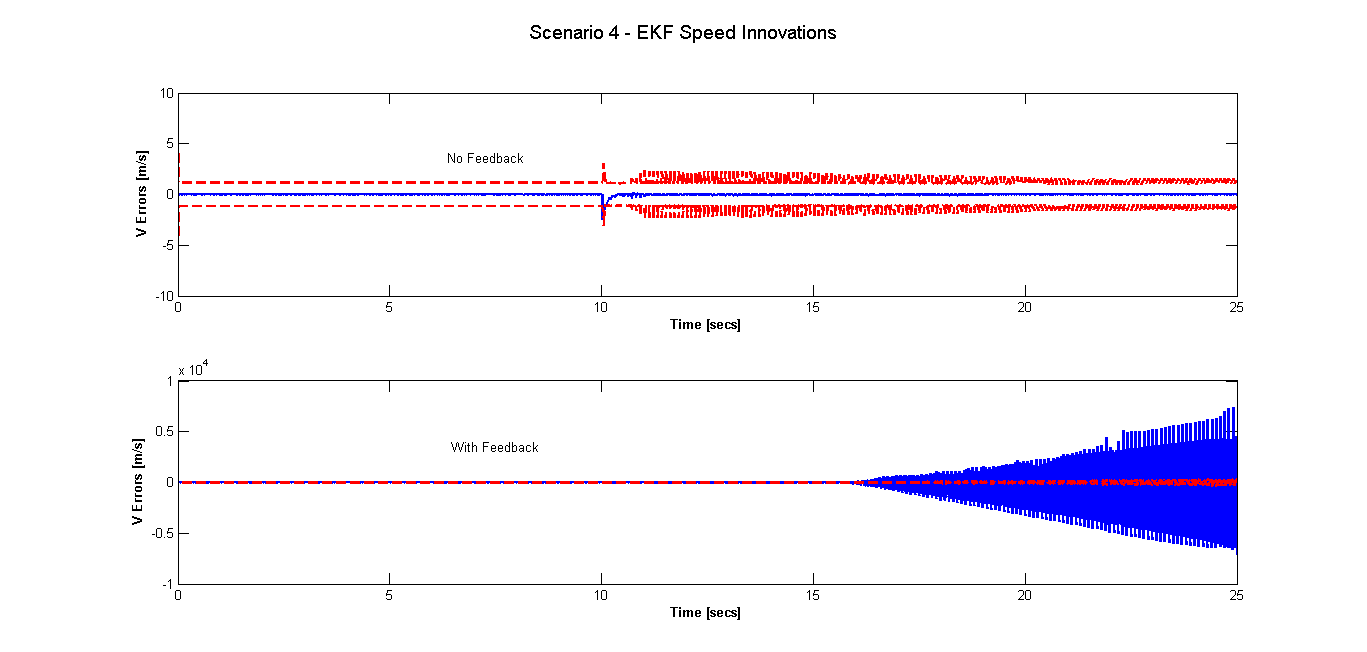}
\caption{Scenario 4 - EKF Velocity Innovations, $2\sigma$ Uncertainty Bounds (red), Velocity Innovations (blue)} 
\label{fig:chap4_scenario4_EKF_velocityErrors}
\end{figure}

\begin{figure}[H]
\center
\includegraphics[scale=0.3]{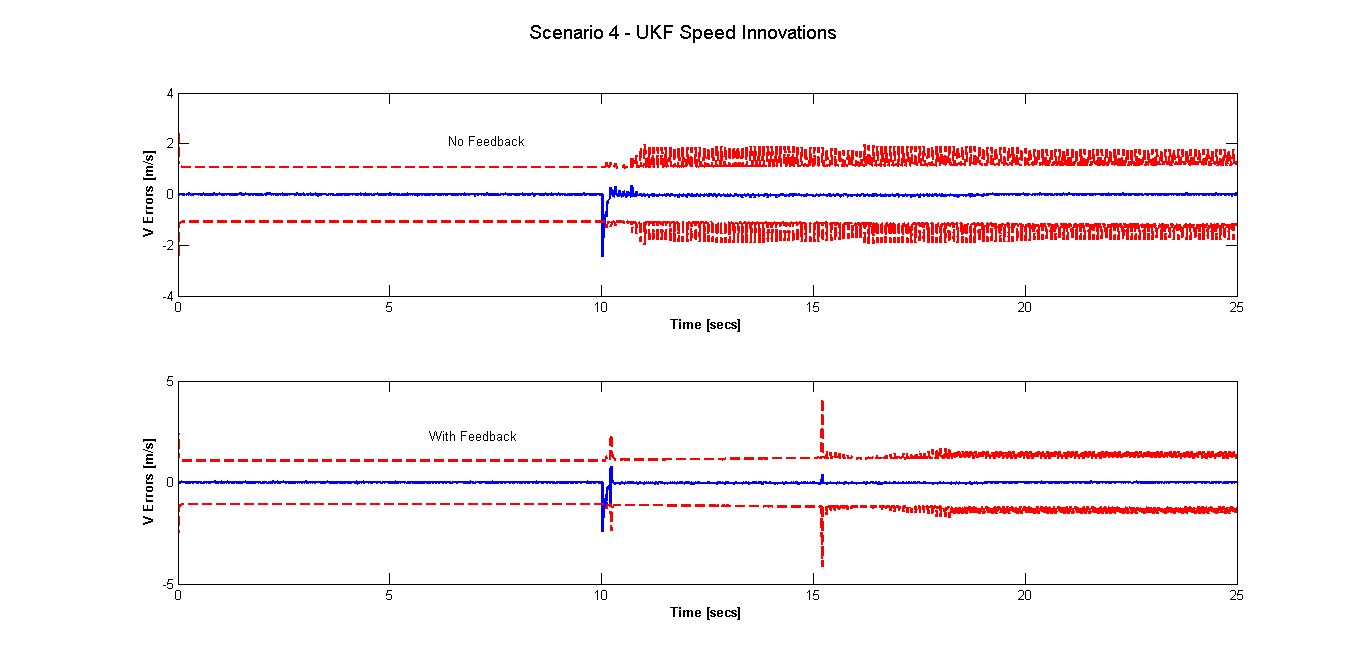}
\caption{Scenario 4 - UKF Velocity Innovations, $2\sigma$ Uncertainty Bounds (red), Velocity Innovations (blue)} 
\label{fig:chap4_scenario4_UKF_velocityErrors}
\end{figure}

\begin{figure}[H]
\center
\includegraphics[scale=0.3]{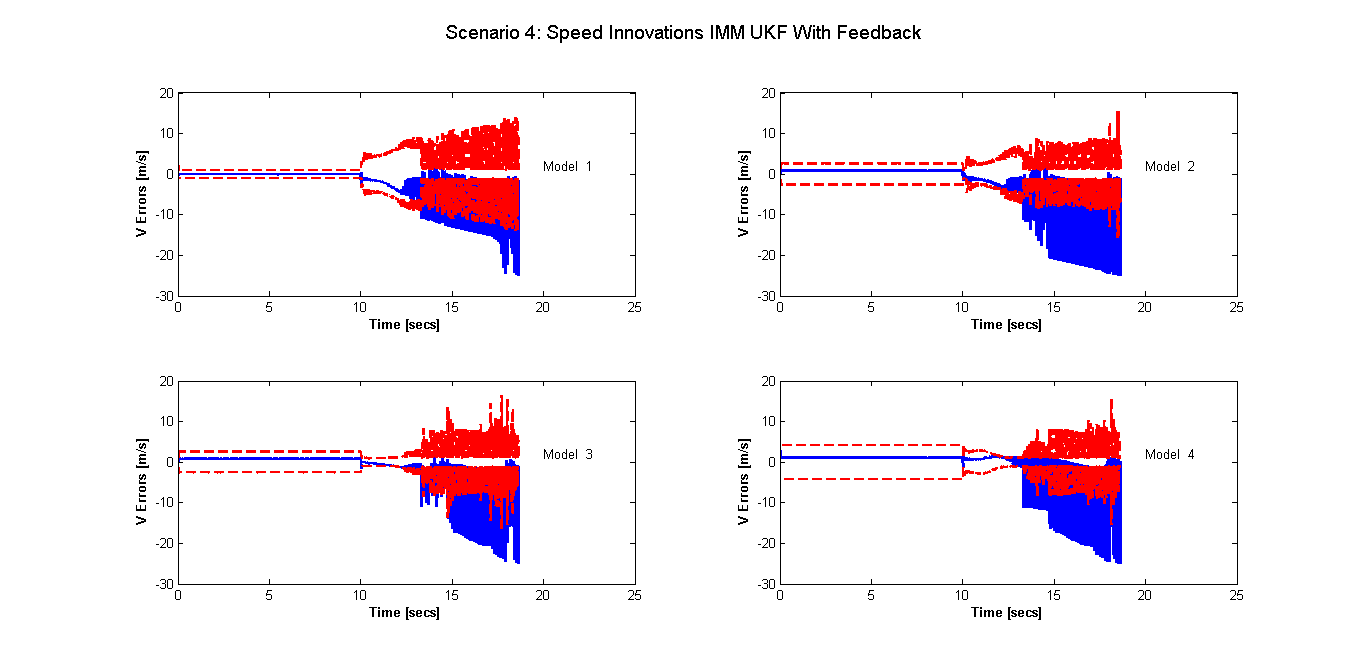}
\caption{Scenario 4 - IMM UKF Velocity Innovations With Feedback, $2\sigma$ Uncertainty Bounds (red), Velocity Innovations (blue)} 
\label{fig:chap4_scenario4_IMMUKF_velocityErrors_WFB}
\end{figure}

Scenario 4 velocity innovations are presented in figures \ref{fig:chap4_scenario4_EKF_velocityErrors} and \ref{fig:chap4_scenario4_UKF_velocityErrors} for the EKF and UKF respectively.  The EKF breaks down 5 seconds after the fault occurs, when feedback is turned on as can be seen by the innovations exceeding the covariance bounds.  Without feedback however the innovations remain within the uncertainty bounds.  The results from the UKF are much better as it produces innovations which remain well within the uncertainty bounds with and without feedback.  Both of the IMM filters failed, because the fault type of scenario 4 was not modelled as a part of the filter design, i.e. the hypothesis for this type of failure is not accounted for and so no model exists for this failure (see figure \ref{fig:chap4_scenario4_IMMUKF_velocityErrors_WFB}).


The wheel radius estimates showed that without a hypothesis on the IMM filters the UKF was the only filter able to produce the correct estimates of the wheel radii.  The UKF results clearly indicated that reconfiguring the controller resulted in a faster convergence to the correct estimate.\\

Plots of wheel speed produced by all four filters with and without feedback show that with this type of fault both wheels were required to work at their constraints.\\ 

The robot trajectories as a result of the different filter information were examined and it was found that none of the filters show full compliance with the reference trajectory.  However the UKF was able to maintain the robot on the path the longest.    

\subsubsection{Filter Re-Design} 
\begin{figure}[H]
\centering
\includegraphics[scale=0.3]{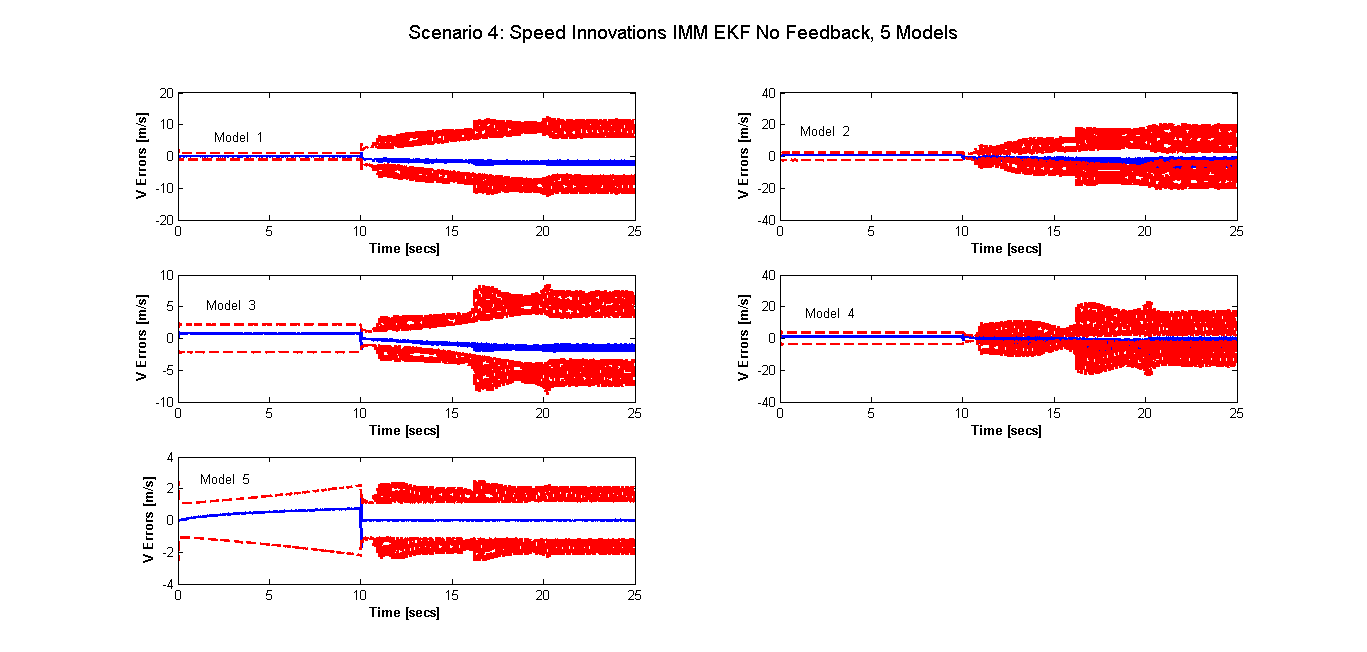}
\caption{Scenario 4 - 5 Model IMM EKF Velocity Innovations No Feedback, $2\sigma$ Uncertainty Bounds (red), Velocity Innovations (blue)} 
\label{fig:chap4_scenario4_IMMEKF_velocityErrors_NFB_5Models}
\end{figure}

\begin{figure}[H]
\centering
\includegraphics[scale=0.3]{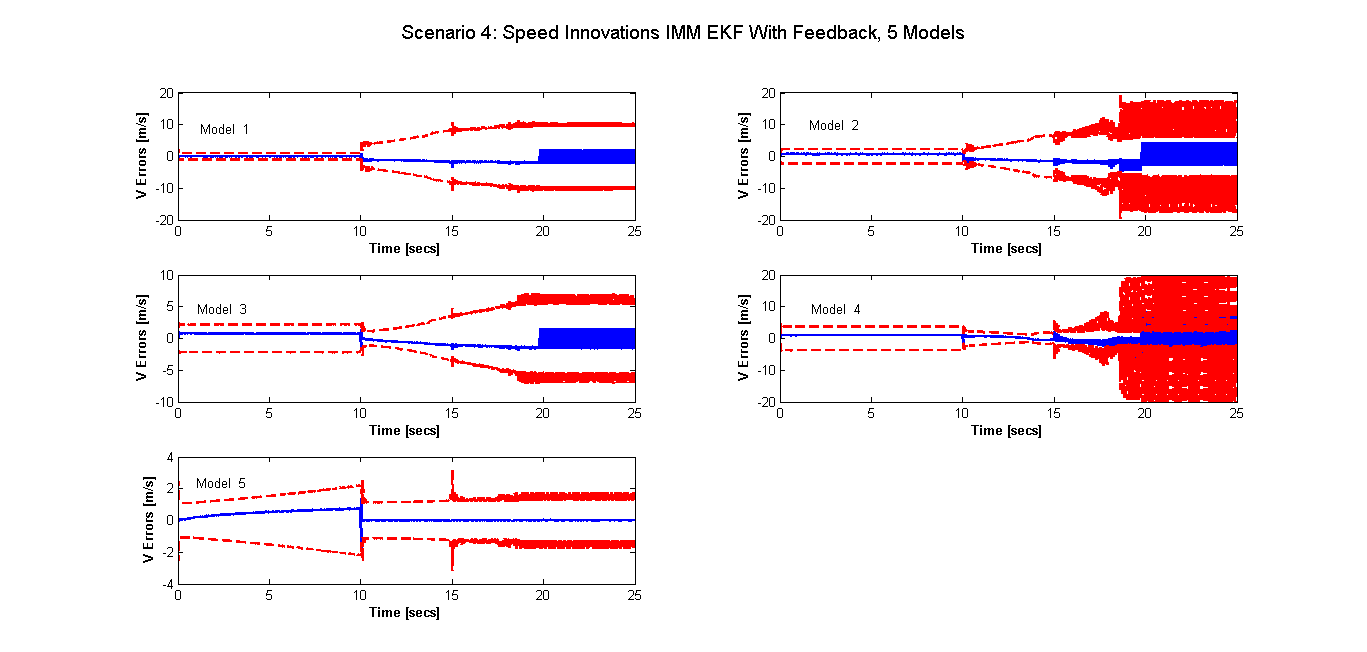}
\caption{Scenario 4 - 5 Model IMM EKF Velocity Innovations With Feedback, $2\sigma$ Uncertainty Bounds (red), Velocity Innovations (blue)} 
\label{fig:chap4_scenario4_IMMEKF_velocityErrors_WFB_5Models}
\end{figure}

- This section looks at the behaviour of the IMM filters by re-designing the filters to accommodate the fault type covered by scenario 4.  The filters were modified by adding a fifth model, which hypothesises the left wheel puncturing in the manner described by scenario 4.  The $Q$ and $R$ matrices remain the same and the initial state vectors for the fifth filter and mode are:
\begin{eqnarray}
\mathbf{x}_5(0) &=& \left[x_0,\,\,\,y_0\,\,\,\psi_0\,\,\,2\,\,\,2\right]^\intercal.
\end{eqnarray}
The mixing probabilities or mode probabilities become:
\begin{equation}
\mu = \left[1/5,1/5,1/5,1/5,1/5\right]^\intercal,
\end{equation}
and the mode transition probabilities matrix $p$ is redefined as:
\begin{equation}
p = \begin{bmatrix}
0.96 & 0.01 & 0.01 & 0.01 & 0.01\\
0.01 & 0.96 & 0.01 & 0.01 & 0.01\\
0.01 & 0.01 & 0.96 & 0.01 & 0.01\\
0.01 & 0.01 & 0.01 & 0.96 & 0.01\\
0.01 & 0.01 & 0.01 & 0.01 & 0.96
\end{bmatrix}.
\end{equation}
The speed innovations for the IMM EKF are presented in figures \ref{fig:chap4_scenario4_IMMEKF_velocityErrors_NFB_5Models} and \ref{fig:chap4_scenario4_IMMEKF_velocityErrors_WFB_5Models} with no feedback and with feedback respectively.  As predicted the results show that if a hypothesis is made the IMM performs very well.  The EKF as part of IMM is able to predict this type of error which it was unable to do as a single filter.  The IMM UKF speed innovation plots showed a higher level of confidence in its estimates compared to its EKF counterpart as the uncertainty is lower and consistent.  Providing feedback in both cases (EKF and UKF IMMs) was shown to increase confidence in the filter estimates.

\begin{figure}[H]
\centering
\includegraphics[scale=0.3]{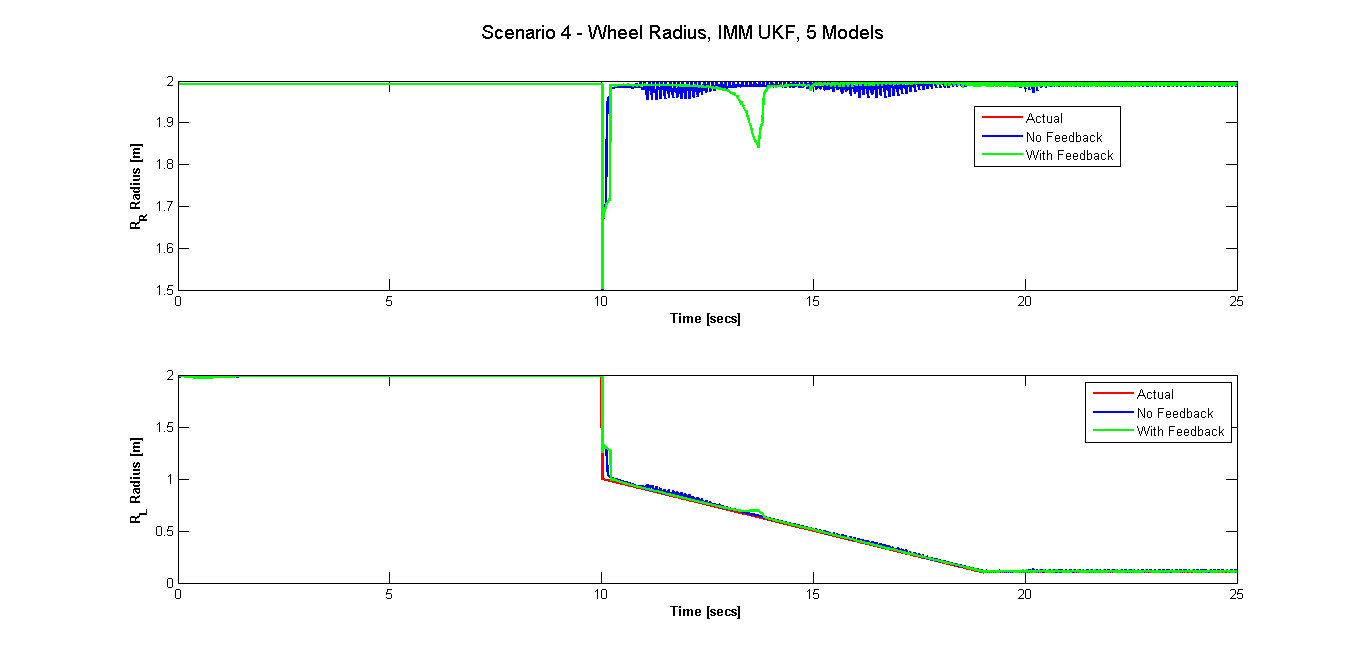}
\caption{Scenario 4 - 5 Model IMM UKF Radius Estimates, Right wheel (top), Left wheel (Bottom)} 
\label{fig:chap4_scenario4_IMMUKF_wheelRadius_5Models}
\end{figure}

The UKF based IMM estimates of the wheel radii for the 5 mode IMM are shown in figure \ref{fig:chap4_scenario4_IMMUKF_wheelRadius_5Models}.  The results of both EKF and UKF based IMMs showed that in both cases the filters do an excellent job of making the correct estimations on the radius of the wheel.
The angular rate plots of the 5 mode IMMs again showed that the control inputs are required to work at the constraints the majority of the time once the fault occurred.  This result is not dependent on the filter type but rather the fault that has occurred makes it impossible for the robot to achieve the desired task while at the same time respecting its constraints.\\   

The trajectory plots showed that for this type of fault where the wheel radius had almost approached zero, the wheels are unable to maintain the reference.  The UKF IMM is again able to keep the robot on the path for a longer time than the EKF IMM.

\subsection{Filter Comparisons}\label{subsection:chap4_FandC}  
\begin{figure}[H]
\centering
\includegraphics[scale=0.3]{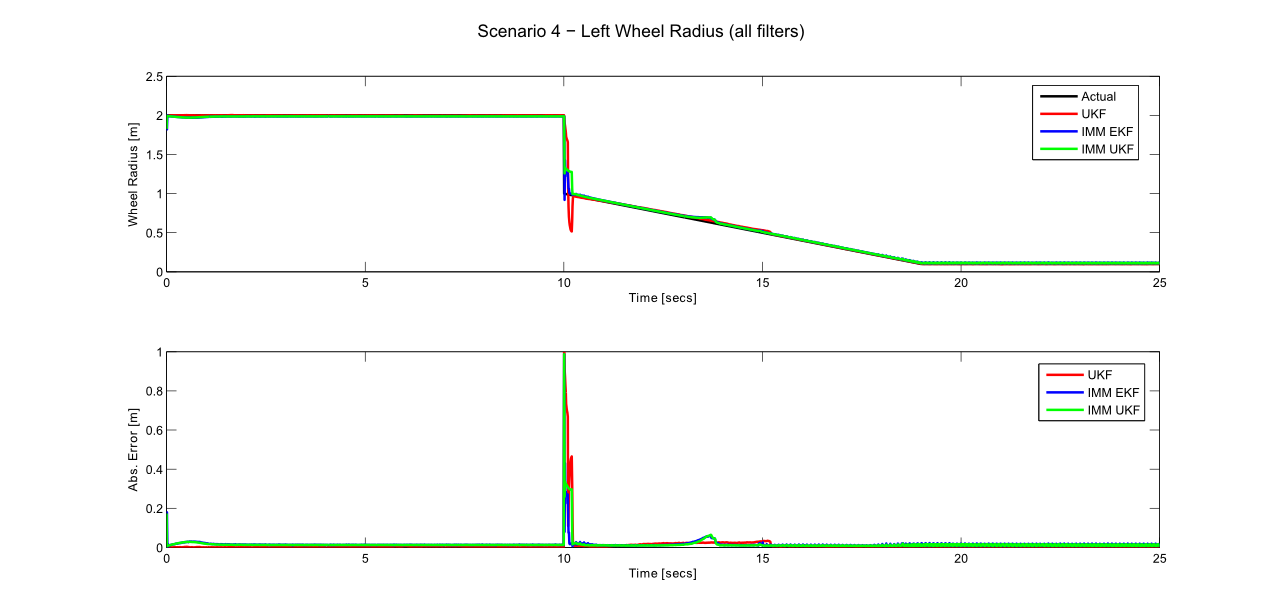} 
\caption{Filter Comparison of Wheel Radius Estimates} 
\label{fig:radiiComparisons}
\end{figure}

The results of the previous section reveal that each of the filters considered exhibits excellent qualities for fault detection and identification.  Figure \ref{fig:radiiComparisons} shows a comparison of the different filters in terms of wheel radius estimation.  The plot shows the radius estimates for scenario 4 which is the worst case scenario out of the 4 scenarios considered.  The output of the EKF and the 4 mode IMM filters have been omitted from the plots as the errors were very high and as a consequence the remaining results were invisible.  The plot shows that in terms of wheel radius estimate performance the UKF performed equally as well as the 5 mode IMM filters.  The IMM filters were able to reach the correct estimate faster than the single UKF particularly when a sudden change is present as is evident at the 10 second mark.  The EKF performance also dramatically increases when used within an IMM configuration.  The drawback of the IMM filter is that all possible scenarios must be accounted for.  The method used in this research illustrated the concept of the IMM, by showing easy adaptability to different situations, and the speed with which it is able to identify and reach the correct estimate.  However predicting exactly how a fault will occur (in this case how a tyre will puncture) is impractical.  A more practical implementation would have been to develop a number of filters each with different process noises that could adapt to all different situations.  The number of filters and the process noise values would have to be determined by trial and error.  In any case an IMM performs well if and only if it is equipped to make a hypothesis on the current situation.  If the given situation is unaccounted for the filter breaks down.  In terms of the single filter, in general the UKF displayed better performance than the EKF, especially in the case of scenario 4 where the non-linearities of the fault caused the single EKF to breakdown.  For these reasons the UKF has been chosen for the final FTC system design.

\subsection{Comparison to Linear MPC}\label{subsec:LMPC}
As a result of the findings given in subsection \ref{subsection:chap4_FandC} only the UKF with feedback is implemented to compare nonlinear MPC with linear MPC.  The results given in the next two sections are for scenarios 2 and 4 respectively however plots for only scenario 2 are given due to space constraints. 

\subsubsection{Scenario 2}
\begin{figure}[H]
\centering
\includegraphics[scale=0.3]{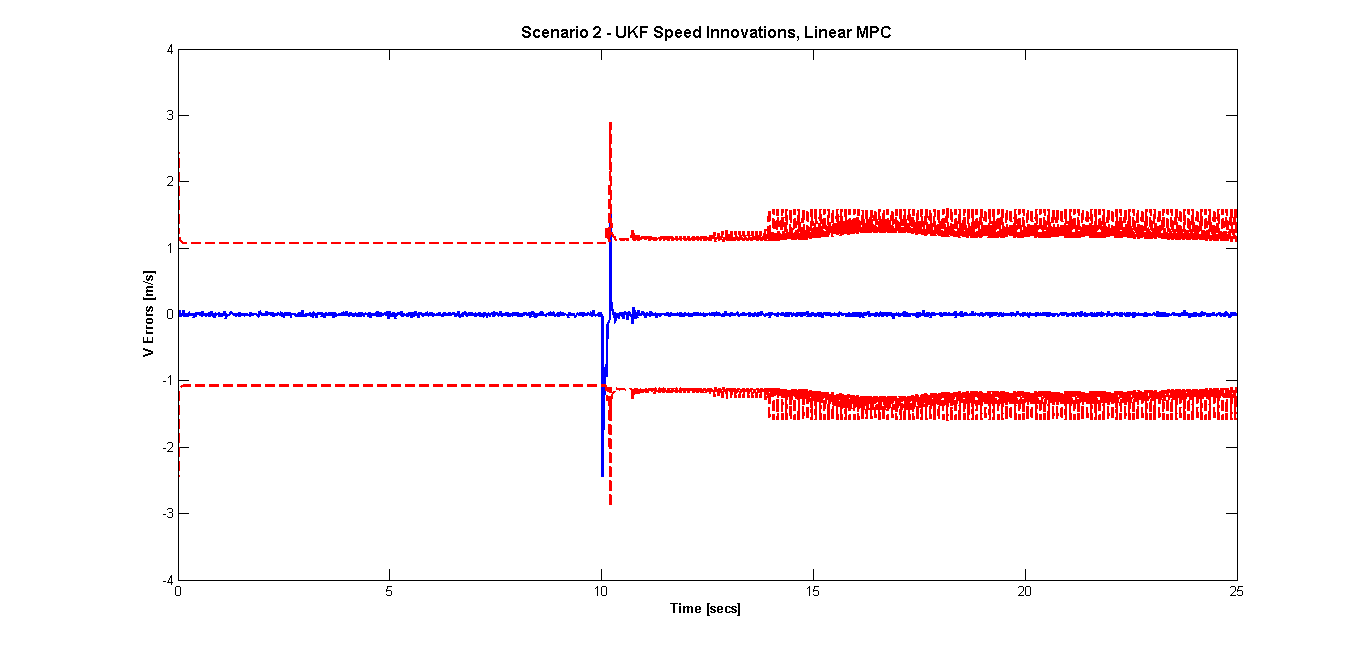}
\caption{Scenario 2 - Linear MPC UKF Speed Innovations, $2\sigma$ Uncertainty Bounds (red), Speed Innovations (blue)} 
\label{fig:chap4_scenario2_UKF_velocityErrors_linearMPC}
\end{figure}

The velocity innovation plots in figure \ref{fig:chap4_scenario2_UKF_velocityErrors_linearMPC} produced by the linear MPC controller show that the innovations remain well within the uncertainty bounds and are approximately zero.  However the uncertainty was seen to be double that produced by the nonlinear controller.  

\begin{figure}[H]
\centering
\includegraphics[scale=0.3]{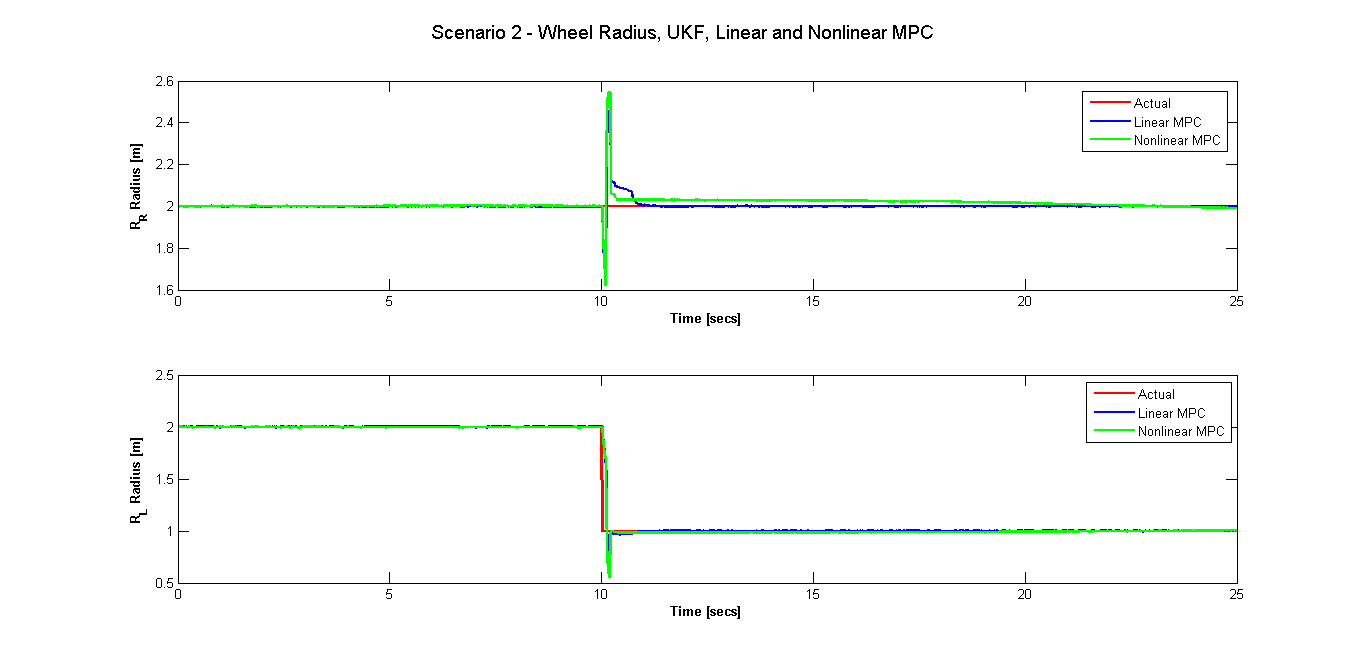}
\caption{Scenario 2 - Comparison of Linear and Nonlinear MPC UKF Radius Estimates, Right wheel (top), Left wheel (Bottom)} 
\label{fig:chap4_scenario2_UKF_wheelRadius_linearMPC}
\end{figure}

The wheel radii plots give in figure \ref{fig:chap4_scenario2_UKF_wheelRadius_linearMPC} show that the estimations produced by a nonlinear controller are the same as those produced by the linear controller. Hence the filter performs well even with linear MPC.

\begin{figure}[H]
\centering
\includegraphics[scale=0.3]{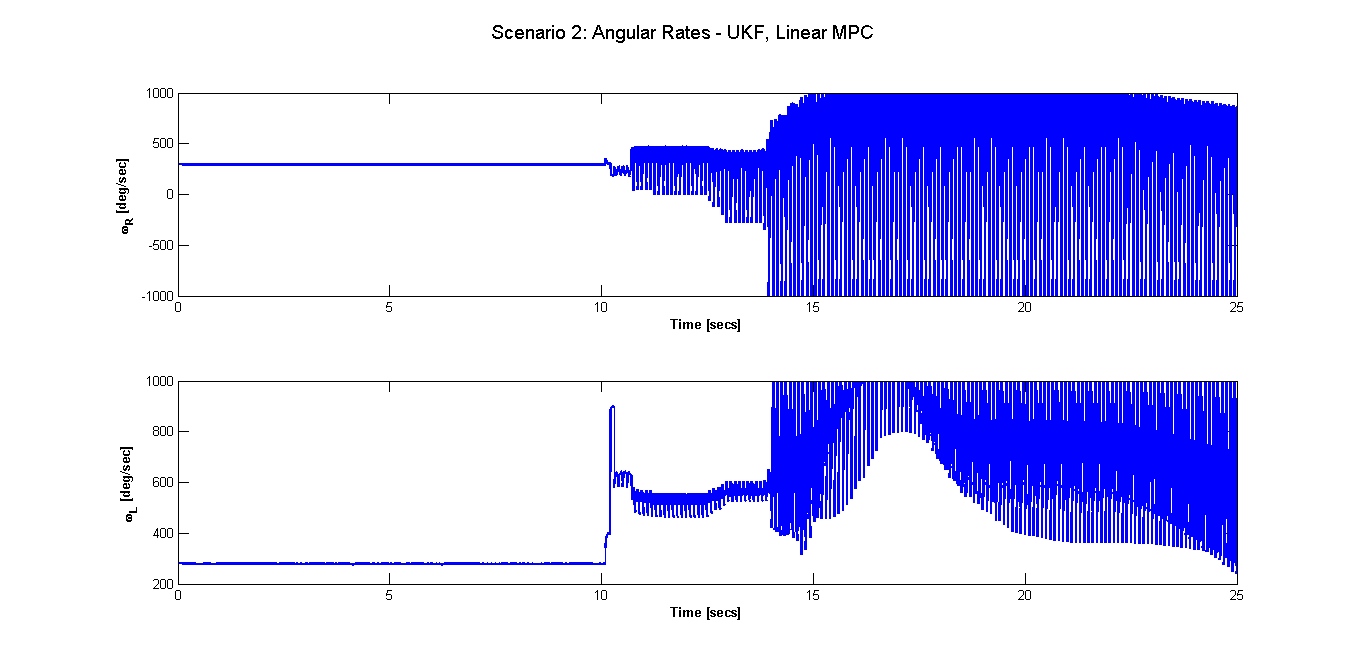} 
\caption{Scenario 2 - Linear MPC UKF Angular Rates, Right wheel (top), Left wheel (Bottom)} 
\label{fig:chap4_scenario2_AR_UKF_linearMPC}
\end{figure}

\begin{figure}[H]
\centering
\includegraphics[scale=0.3]{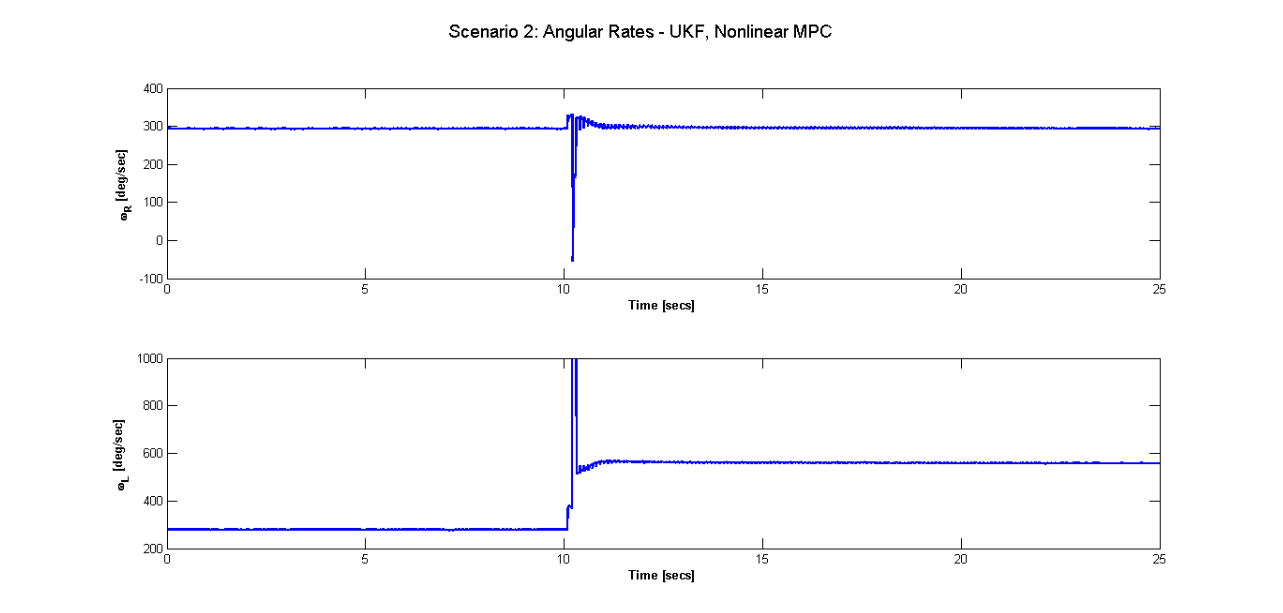} 
\caption{Scenario 2 - Nonlinear MPC UKF Angular Rates, Right wheel (top), Left wheel (Bottom)} 
\label{fig:chap4_scenario2_AR_UKF_NMPC}
\end{figure}

\begin{figure}[H]
\centering
\includegraphics[scale=0.3]{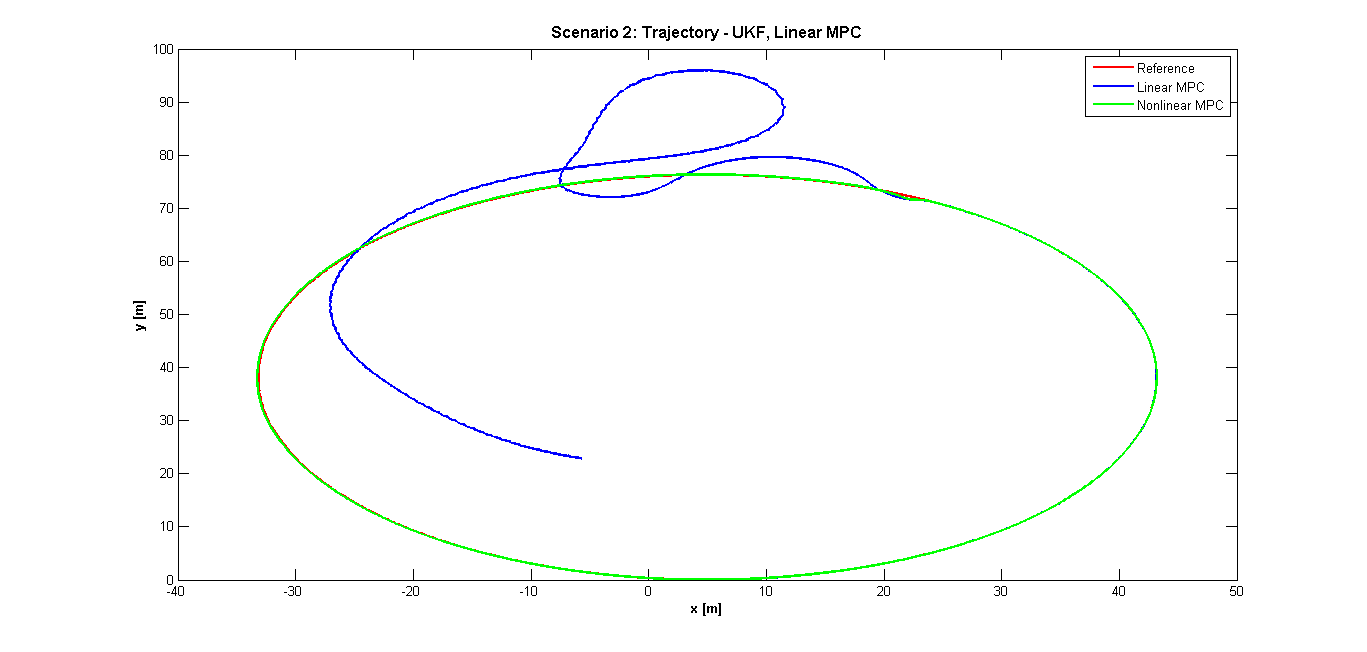} 
\caption{Scenario 2 - Linear and Nonlinear MPC UKF Trajectory} 
\label{fig:chap4_scenario2_traj_UKF_linearMPC}
\end{figure}

The angular rate plots for the linear MPC controller (figure \ref{fig:chap4_scenario2_AR_UKF_linearMPC}) show that five seconds after the fault occurred the linear controller pushes the wheels to operate at their constraints and is unable to tolerate the faulty condition.  The angular rates produced by the NMPC controller (figure \ref{fig:chap4_scenario2_AR_UKF_NMPC}) show that it was able to easily adapt to the fault.  This is further illustrated in the trajectory plot given in figure \ref{fig:chap4_scenario2_traj_UKF_linearMPC} which clearly shows that the nonlinear MPC controller does an excellent job of keeping the robot on the path despite a faulty wheel whereas the solution produced by the linear controller has diverged.   

A point to note here is that while the trajectory tracking is good, the switching behaviour observed in wheel rotation rates is highly undesirable, and in the real world would produce wheel slippage, and high levels of wear on tyres and mechanicals. The wheel rates produced by the nonlinear MPC controller show a dramatic decrease in this switching behaviour.  It is, however, still present (figure \ref{fig:chap4_scenario2_AR_UKF_NMPC}).  This limit cycle behaviour is a negative characteristic that needs to be addressed before this technique can be applied to real systems.  This exercise was purely for proof of concept and is not a practical application.  The model is entirely kinematic and for it to represent a more practical scenario further work is required to eliminate the limit cycling behaviour; for example adding actuation activity to the cost function and an angular acceleration term to avoid wheel slippage issues.

\subsubsection{Scenario 4}
The speed innovations were plotted for scenario 4 and showed similar trends to those above in that the innovations are quite small; however the uncertainties with the linear controller are higher than those produced as a result of nonlinear MPC. The estimates of the wheel radii however are very good and were seen to be the same for both linear and nonlinear controllers.\\

As expected the linear controller was unable to maintain the robot on the path.  Neither of the controllers were able to drive the robot to the path as the wheel radius had almost reached 0m making it infeasible for the robot to continue.  The angular rates showed that the linear controller constantly demanded operation at the constraints oscillating between the upper and lower limits continuously.  The nonlinear MPC controller results for wheel angular rates showed that the right wheel oscillates between the upper and lower bounds, however the left wheel is required to constantly work at the upper bound.

\section{Conclusion}\label{section:chap4_conclusion}
The analysis from this work has proven the feasibility of the NMPC controller design with filter estimates for controller reconfiguration as a viable solution to fault tolerant control.  Four different filters were compared, the EKF, the UKF, the IMM EKF and the IMM UKF filters.  The results showed that in terms of fault detection performance the UKF is the best candidate in a trajectory tracking scenario.\\

Comparisons were also made between the performance of nonlinear MPC and linear MPC.  The results clearly show that for the purposes of reconfigurable fault tolerant control the nonlinear MPC controller has better performance.\\

The next phase of this research is to implement the NMPC pseudospectral controller with a UKF based FDI subsystem to an aircraft, for fault tolerant flight control.


\end{document}